\documentclass[10pt]{article}
\usepackage{amssymb, amsmath, url, graphicx, setspace, geometry}
\usepackage{calrsfs}
\usepackage{wasysym, xcolor}
\def\3{\subset }
\def\4{\subseteq }
\def\<{\left<}
\def\>{\right>}

\def\bit{\begin{itemize}}
\def\eit{\end{itemize}}
\def\3{\subset }
\def\4{\subseteq }

\def\0{\leqno}

\def\barr{\begin{array}}
\def\earr{\end{array}}

\def\Z{{\rlap{$\kern2pt{\rm Z}$}{\rm Z}\,}}

\title{\bf A set theoretic version of equations on groups}
\author{Mihai-Silviu Lazorec}
\date{February 28, 2026}

\begin{document}

\maketitle

\begin{abstract}
Let $G$ be a finite group. The aim of this paper is to study the number of solutions $S\subseteq G$ of the equation $\mho^{\{n\}}(S)=L$, where $L$ is a non-empty subset of $G$, $n$ is a positive integer and $\mho^{\{n\}}(S)=\{ s^n \ | \ s\in S\}$. Besides our  findings obtained in this general frame, we also outline some results which hold for some particular cases such as: \textit{i)} $L$ is a normal subset of $G$; \textit{ii)} $G$ is abelian; \textit{iii)} $G$ is an extraspecial $p$-group.
\end{abstract}



\section{Introduction}

Let $G$ be a finite group, $n$ be a positive integer and $l$ be an element of $G$. We denote by $N_{n,l}$ the set of solutions $x\in G$ of the equation $x^n=l$. One of the best-known results concerning such an equation is due to Frobenius (see \cite{2, 4}):\\

\textbf{Theorem 1.1.} \textit{Let $G$ be a finite group, $n$ be a positive integer and $l\in G$. Then $|N_{n,l}|$ is a multiple of $(n, |C_G(l)|)$.} \\

While the  most popular version of Theorem 1.1 is the one corresponding to the case $l=1$ (check \cite{11} for a proof and \cite{9} for some applications), the literature also includes various generalizations of this result (see \cite{6}). We recall the following two results which outline what happens if $l$ ranges over a subgroup or a normal subset of $G$.\\

\textbf{Theorem 1.2.} \textit{Let $G$ be a finite group and $n$ be a positive integer.
\begin{itemize}
\item[i)] If $L$ is a subgroup of $G$, then $\sum\limits_{l\in L}|N_{n, l}|$ is a multiple of $|L|$ (see Theorem 5 of \cite{12} or Corollary 4 of \cite{10});
\item[ii)] If $L$ is a normal subset of $G$, then $\sum\limits_{l\in L}|N_{n, l}|$ is a multiple of $(n, |G|)$ (see the first two  statements in section 5, p. 163 of \cite{18} ).
\end{itemize}}

An equivalent perspective on the previous theorem is that it provides some information on the number of elements $x\in G$ such that $x^n\in L$ when $L$ is a subgroup or a normal subset of $G$. In general, for a normal subgroup $L$ of $G$, one can not state that one of the two items of Theorem 1.2 constitutes a stronger result than the other one since $|L|$ and $(n, |G|)$ are not comparable (with respect to the usual divisibility relation). Indeed, if $G\cong S_3, L\cong\langle (1, 2, 3)\rangle$ and $n=2$, then $\sum\limits_{l\in L}|N_{n, l}|=6$ is a multiple of both $|L|=3$ and $(n,|G|)=2$.

Let $G$ be a finite group and $n$ be a positive integer. We denote by $\mathcal{P}^*(A)$ the set formed of the non-empty subsets of a set $A$. For $S\in \mathcal{P}^*(G)$, we denote by $\mho^{\{ n\}}(S)$ the set containing the $n$th powers of the elements of $S$, i.e. $\mho^{\{ n\}}(S)=\{ s^n \ | \ s\in S\}$. For $L\in \mathcal{P}^*(G)$, we are interested in obtaining some results on the size of the set $N_{n,L}$ of solutions $S$ of the equation
$$\mho^{\{ n\}}(S)=L.$$
In this general case, we show that $|N_{n,L}|$ is a multiple of $\big[2^{(n, |C_G(L)|)}-1\big]^{|L|}$. Variations of this result are obtained when $L$ is a normal subset of $G$ or $f:G\longrightarrow G$ given by $f(x)=x^n, \ \forall \ x\in G,$ is a group homomorphism. Also, if $G$ is abelian or an extraspecial $p$-group, we determine the exact value of $|N_{n,L}|$ and we check if it divides $|\mathcal{P}^*(G)|^{|L|}$ (if $|N_{n,L}|\neq 0$).

We end this section by enumerating some common notations that are going to be used in this paper. $D_8$ stands for the dihedral group with 8 elements, $Q_8$ denotes the quaternion group and $C_k$ is the cyclic group with $k\geq 1$ elements. The group
$$M_{p^k}=\langle x, y \mid x^{p^{k-1}}=y^p=1, y^{-1}xy=x^{p^{k-2}+1}\rangle$$ is a  modular $p$-group of order $p^k$ ($k\geq 3$ if $p$ is odd; $k\geq 4$ if $p=2$). If $p$ is odd, the non-abelian $p$-group of order $p^3$ and exponent $p$ is
$$He_p=\langle x, y, z \mid x^p=y^p=z^p=1, [x, z]=[y, z]=1, [x, y]=z\rangle.$$ If $G$ and $H$ are finite groups, we denote their central product with joint centers by $G\circ H$, while $G^{\circ k}$ stands for the central product of $k$ copies of $G$ (if $k=0$, then $G^{\circ k}=\{1\}$; if $k=1$, then $G^{\circ k}=G$). The order of an element $x\in G$ is denoted by $o(x)$, $n_d(G)$ is the number of elements of order $d$ in $G$ and $exp(G)$ is the exponent of $G$.  

\section{Main results}

Let $G$ be a finite group, $n$ be a positive integer, $L\in \mathcal{P}^*(G)$ and $l$ be an element of $L$. As a preliminary result, we precisely indicate $|N_{n, L}|$ when some arithmetic properties referring to $n$ and $exp(G)$ hold. We also show that $|N_{n,l}|$ is a multiple of a quantity which is independent of $l$. This is achieved by weakening (with respect to divisibility) the g.c.d. that appears in Theorem 1.1. Finally, we establish a connection between the function $f$ introduced in section 1 and $N_{n, L}$.\\

\textbf{Lemma 2.1.} \textit{Let $G$ be a finite group, $n$ be a positive integer, $L\in \mathcal{P}^*(G)$ and $l\in L$. The following statements hold: 
\begin{itemize}
\item[i)] If $exp(G)\mid n$, then 
$
|N_{n, L}|=\begin{cases} 0 & \text{, if } L\neq\{ 1\} \\ 2^{|G|}-1 & \text{, if } L=\{1 \}
\end{cases}.
$
\item[ii)] If $(n, exp(G))=1$, then $|N_{n, L}|=1$;
\item[iii)] $|N_{n,l}|$ is a multiple of $(n, |C_G(L)|)$;
\item[iv)] Assume that the function $f:G\longrightarrow G$  given by $f(x)=x^n, \ \forall  \ x\in G$, is a group homomorphism. If $L\not\subseteq Im(f)$, then $N_{n, L}=\emptyset$; otherwise, $|N_{n, L}|=[2^{|Ker(f)|}-1]^{|L|}.$
\end{itemize}}

\textbf{Proof.} The proof of \textit{i)} is immediate. Hence, we only share some details on the proofs of the remaining items. 
 
\textit{ii)} Under the given hypothesis, there is a positive integer $m$ such that $nm \equiv  1 \ (mod \ exp(G))$. It is easy to check that $\mho^{\{m\}}(L)\in N_{n, L}$.
For a given positive integer $a$, we consider the function 
$$F_a:\mathcal{P}^*(G)\longrightarrow \mathcal{P}^*(G) \text{ given by } F_a(S)=\mho^{\{a\}}(S), \ \forall \ S\in\mathcal{P}^*(G).$$
We easily infer that $F_n$ is a bijection  whose inverse is $F_m$. Hence, the solution $\mho^{\{m\}}(L)$ is unique and the conclusion follows.   

\textit{iii)} Since $C_G(L)\leq C_G(l)$, we deduce that $|C_G(L)|\mid |C_G(l)|$. Then $(n, |C_G(L)|)\mid (n, |C_G(l)|)$ and the conclusion follows as a consequence of the fact that $(n, |C_G(l)|)\mid |N_{n,l}|$ by Theorem 1.1. 

\textit{iv)} The first part of the statement is clear. Let $L=\{ l_1, l_2, \ldots, l_k\}\subseteq Im(f)$, where $k$ is a positive integer. Then, for any $i\in \{ 1, 2, \ldots, k\}$, we can choose $x_i\in N_{n, l_i}$. It is easy to check that $f^{-1}(L)=\bigcup\limits_{i=1}^{k}x_iKer(f)$. Then
$$S\in N_{n, L}\Longleftrightarrow S=\bigcup\limits_{i=1}^k S_i, \text{ where } S_i\in\mathcal{P}^*(x_iKer(f)), \ \forall \ i\in\{ 1, 2, \ldots, k\}.$$
Since the sets $S_i$ are disjoint and each of them can be chosen in $2^{|Ker(f)|}-1$ ways, we obtain the stated result.
\hfill\rule{1,5mm}{1,5mm}\\ 

As a consequence of Lemma 2.1 \textit{iii)}, it is clear that $\sum\limits_{l\in L}|N_{n,l}|$ is a multiple of $(n, |C_G(L)|)$, for any  $L\in \mathcal{P}^*(G)$. It is known that $C_G(L)\leq N_G(L)$. For any normal subset $L\in \mathcal{P}^*(G)$, since $N_G(L)=G$, we know that $\sum\limits_{l\in L}|N_{n,l}|$ is a multiple of $(n, |N_G(L)|)$ according to Theorem 1.2  \textit{ii)}. We pose the following question related to these aspects.\\

\textbf{Open problem.} \textit{Let $G$ be a finite group, $n$ be a positive integer and $L$ be a non-normal subgroup of $G$. Is it true that $\sum\limits_{l\in L}|N_{n,l}|$ is a multiple of $(n, |N_G(L)|)$?}\\

A check performed via GAP \cite{20} shows that the answer is positive for any group up to order 255. Related to this question, by Theorem 1.2 \textit{ii)},  we note that the number of elements $x\in N_G(L)$ such that $x^n\in L$ is a multiple of $(n, |N_G(L)|)$. However, this result does not hold if we search for solutions $x$ in a coset $gN_G(L)\neq N_G(L)$. For instance, if we take $G\cong S_3, L\cong \langle (2,3) \rangle, n=2$ and $g=(1, 3)$, then $(n, |N_G(L)|)=2$ and we find one solution $x\in gN_G(L)$ such that $x^n\in L$. Hence, by considering that the answer to the open question is affirmative, this can not be justified by showing that each coset $gN_G(L)$ contains a number of solutions which is a multiple of $(n, |N_G(L)|)$. This would have indeed implied that $\sum\limits_{l\in L}|N_{n,l}|$ is a multiple of $(n, |N_G(L)|)$. 

A recent perspective on the number of elements $x\in G$ such that $x^n\in L$, where $L$ is a subgroup of $G$, was outlined in \cite{17}. More exactly, $$x^n\in L\Longleftrightarrow o_L(x)\mid n,$$ where $o_L(x)$ is the order of $x$ relative to $L$, i.e. $o_L(x)=\min\{k\in\mathbb{N}^* \mid  x^k\in L\}$. Other results concerning this generalization of the classic concept of element order are given in \cite{14, 19}. By using the above equivalence, one could express the sum $\sum\limits_{l\in L}|N_{n, l}|$ as follows.\\

\textbf{Proposition 2.2.} \textit{Let $G$ be a finite group, $n$ be a positive integer and $L$ be a subgroup of $G$. Then $$\sum\limits_{l\in L}|N_{n, l}|=|\{x\in G \mid  o_L(x)\mid n\}|.$$}

We may evaluate the number of solutions $S\subseteq G$ of the equation $\mho^{\{n\}}(S)=L$, where $L\in\mathcal{P}^*(G)$.\\ 

\textbf{Theorem 2.3.} \textit{Let $G$ be a finite group, $n$ be a positive integer and $L\in \mathcal{P}^*(G)$. Then $|N_{n,L}|$ is a multiple of $\big[2^{(n, |C_G(L)|)}-1\big]^{|L|}$. In particular, if $|N_{n, L}|\neq 0$, then $|N_{n, L}|$ is an odd number.}\\

\textbf{Proof.} If there exists $l\in L$ such that $N_{n,l}=\emptyset$, then $N_{n,L}=\emptyset$, so the conclusion holds.

Suppose that $N_{n,l}\neq\emptyset$, for all $l\in L$. It is easy to observe that 
$$S\in N_{n, L} \Longleftrightarrow S=\bigcup\limits_{l\in L}S_l, \text{ where } S_l\in\mathcal{P}^*(N_{n, l}), \ \forall \ l\in L.$$ Let $l\in L$. Then $S_l$ can be chosen in $2^{|N_{n,l}|}-1$ ways. By Lemma 2.1 \textit{iii)}, we know that $(n, |C_G(L)|)$ divides $|N_{n, l}|$. Consequently $2^{(n, |C_G(L)|)}-1\mid 2^{|N_{n,l}|}-1$, so the number of ways in which $S_l$ can be chosen is a multiple of $2^{(n, |C_G(L)|)}-1$. Since the sets $S_l$ are disjoint, it follows that $$|N_{n, L}|=\prod\limits_{l\in L}(2^{|N_{n,l}|}-1)$$ 
is a multiple of $[2^{(n, |C_G(L)|)}-1]^{|L|}$ and it is also an odd integer, as desired. 
\hfill\rule{1,5mm}{1,5mm}\\

By using the result highlighted in Theorem 1.2 \textit{ii)}, a similar reasoning as the one done in the previous proof could be repeated to justify the validity of the following corollary.\\

\textbf{Corollary 2.4.} \textit{Let $G$ be a finite group and $n$ be a positive integer. If $L$ is a normal subset of $G$, then $|N_{n, L}|$ is a multiple of $\prod\limits_{i=1}^k\big[2^{(n, |C_G(L_i)|)}-1\big]^{|L_i|}$, where $\{L_1, L_2, \ldots, L_k\}$ is the partition of $L$ into conjugacy classes.}\\

If $L$ is a normal subset of $G$, it is easy to show that $\prod\limits_{i=1}^k\big[2^{(n, |C_G(L_i)|)}-1\big]^{|L_i|}$ is a multiple of $\big[2^{(n, |C_G(L)|)}-1\big]^{|L|}$. Hence, in this case, Corollary 2.4 is a stronger result than Theorem 2.3. Another stronger result holds if the function $f$ introduced in section 1 is a group homomorphism.\\

\textbf{Corollary 2.5.} \textit{Let $G$ be a finite group, $n$ be a positive integer and  $L\in \mathcal{P}^*(G)$. If $f:G\longrightarrow G$ given by $f(x)=x^n, \ \forall \ x \in G$, is a group homomorphism,   then $|N_{n,L}|$ is a multiple of $\big[2^{(n, |G|)}-1\big]^{|L|}$.}\\

\textbf{Proof.} Let $l=1$ in Theorem 1.1. Then $(n, |G|)\mid |Ker(f)|$. Then $[2^{|Ker(f)|}-1]^{|L|}$ is a multiple of $[2^{(n, |G|)}-1]^{|L|}$. The proof is complete by  Lemma 2.1 \textit{iv)}.
\hfill\rule{1,5mm}{1,5mm}\\

If $L$ is not a normal subset of $G$ and $f$ is not a group homomorphism, then the previous two corollaries do not hold. For instance, if $G\cong S_3$, $n=2$ and $L=\{ (), (1, 2, 3)\}$, then the equation $\mho^{\{n \}}(S)=L$ has 15 solutions $S$, while $\big[2^{(n, |G|)}-1\big]^{|L|}=9$. 

In what follows, we aim to determine exact formulas for $|N_{n, L}|$ for some specific classes of groups. We begin by considering the case of finite abelian groups.\\

\textbf{Theorem 2.6.} \textit{Let $G\cong G_1\times G_2\times\ldots\times G_k$ be a finite abelian group, where $k\in\mathbb{N}^*$ and $G_i\cong C_{p_i^{\alpha_{i1}}}\times C_{p_i^{\alpha_{i2}}}\times\ldots \times C_{p_i^{\alpha_{iv_i}}}$ is the Sylow $p_i$-subgroup of $G$ for any $i\in\{1,2,\ldots, k \}$. Let $n$ be a positive integer and $L\in\mathcal{P}^*(G)$. Then $$|N_{n, L}|\in\bigg\{ 0, \bigg[2^{\prod\limits_{i=1}^k\prod\limits_{j=1}^{v_i}(n, p_i^{\alpha_{ij}})}-1\bigg]^{|L|}\bigg\}.$$
In particular, if $|N_{n,L}|\neq 0$, then $|N_{n,L}|\mid |\mathcal{P}^*(G)|^{|L|}.$}\\

\textbf{Proof.} Since $G$ is abelian, the function $f:G\longrightarrow G$ given by $f(x)=x^n, \ \forall \ x\in G$, is a group homomorphism. Hence, by Lemma 2.1 \textit{iv)}, we have $|N_{n, L}|=0$ or it suffices to evaluate $|Ker(f)|$. 

Let $x=(g_{11}^{\beta_{11}}, g_{12}^{\beta_{12}}, \ldots, g_{1v_1}^{\beta_{1v_1}}, g_{21}^{\beta_{21}}, g_{22}^{\beta_{22}}, \ldots, g_{2v_2}^{\beta_{2v_2}}, \ldots,  g_{k1}^{\beta_{k1}}, g_{k2}^{\beta_{k2}}, \ldots, g_{kv_k}^{\beta_{kv_k}})$ be an element of $G$, where $0\leq \beta_{ij}<p^{\alpha_{iv_i}}$ for all $i\in\{1,2,\ldots, k\}, j\in \{1, 2, \ldots, v_i\}$. Then 
$$x\in Ker(f)\Longleftrightarrow \forall \ i\in\{1,2,\ldots, k\}, j\in \{1, 2, \ldots, v_i\}: n\beta_{ij}\equiv 0 \ (mod \ p_{i}^{\alpha_{ij}}).$$
For fixed $i$ and $j$, the previous congruence has $(n, p_i^{\alpha_{ij}})$ solutions. Hence, $|Ker(f)|=\prod\limits_{i=1}^k\prod\limits_{j=1}^{v_i}(n, p_i^{\alpha_{ij}})$ and we are done.

In what concerns the particular conclusion, if $|N_{n, L}|\neq 0$, since $|Ker(f)|\mid |G|$, it follows that $|N_{n, L}|\mid |\mathcal{P}^*(G)|^{|L|}.$
\hfill\rule{1,5mm}{1,5mm}\\ 

In general, if $|N_{n,L}|\neq 0$, then the divisibility $|N_{n,L}|\mid |\mathcal{P}^*(G)|^{|L|}$ does not hold beyond the class of finite abelian groups. For instance, if $G\cong D_8, n=2, L=\lbrace 1\rbrace$, then $|N_{n,L}|=63$ and $|\mathcal{P}^*(G)|^{|L|}=255.$ We pose the following question related to this matter.\\

\textbf{Open problem.} \textit{Let $\mathcal{F}$ be the class of all finite groups. For $G\in\mathcal{F}$, we consider the set $A_G=\{ (n, L)\in \mathbb{N}^*\times \mathcal{P}^*(G) \mid |N_{n,L}|\neq 0\}$. Determine the class 
$$\mathcal{D}=\big\{G\in\mathcal{F} \mid |N_{n, L}|\mid |\mathcal{P}^*(G)|^{|L|}, \ \forall \ (n,L)\in A_G\big\}.$$}

According to Theorem 2.6, any finite abelian group belongs to $\mathcal{D}$. The first non-abelian finite groups which have the same property are $C_2^2\rtimes C_4$ (SmallGroup(16,3)), $C_4\rtimes C_4$ (SmallGroup(16,4)), $M_{16}$ (SmallGroup(16,6)) and $D_8\circ C_4$ (SmallGroup(16,13)).

Related to the same question, let $G$ be a finite $p$-group with $exp(G)=p^k$, where $k\in\mathbb{N}^*$, such that $G\in\mathcal{D}$. For $L=\{ 1\}$ and $n=p^i$, with $i\in\{ 1,2,\ldots, k\}$, we get that
$$S\in N_{n, L}\Longleftrightarrow S\in \mathcal{P}^*(\Omega_{\{i\}}(G)), \text{ where } \Omega_{\{i\}}(G)=\{x\in G \mid o(x)\mid p^i\}.$$ Consequently, we get
$$|N_{n,L}|=2^{|\Omega_{\{i\}}(G)|}-1.$$
Since $|N_{n, L}|\mid |\mathcal{P}^*(G)|^{|L|}$, it follows that $2^{|\Omega_{\{i\}}(G)|}-1\mid 2^{|G|}-1$, so 
$$|\Omega_{\{i\}}(G)|\mid |G|, \ \forall \ i\in\{1,2,\ldots, k \}.$$
Therefore, if a finite $p$-group $G$ with $exp(G)=p^k$ belongs to $\mathcal{D}$, then $|\Omega_{\{i\}}(G)|\mid |G|$, for any $i\in\{1,2,\ldots, k \}$. For instance, such a set of divisibilities holds for any finite regular $p$-group (see section 4.13 of \cite{5}) and any powerful $p$-group, with $p$ being an odd prime in the latter case (see \cite{7, 15, 16}), because $\langle \Omega_{\{i\}}(G)\rangle=\Omega_{\{i\}}(G),$ for all $i\in \{1,2,\ldots, k \}.$ For an odd prime $p$, it is known that any extraspecial $p$-group is regular. This is guaranteed by 10.2 Satz \textit{c)} of \cite{8}. Hence, it is worthwhile studying whether such a group belongs to $\mathcal{D}$. The case $p=2$ is treated as well.

We recall that a finite $p$-group $G$ is called extraspecial if $G$ is of class 2, $G'=Z(G)=\Phi(G)$ and $|G'|=p$ (p.183 of \cite{3}). A description of the extraspecial $p$-groups is given by Theorem 2.3 of \cite{1}. Let $k\geq 1$ be an integer and $G$ be an extraspecial $p$-group. It is known that $|G|=p^{2k+1}$ and
\begin{itemize}
\item[\textit{i)}] if $p$ is odd, then $G\cong M_{p^3}^{\circ k}$ or $G\cong He_p^{\circ k}$;
\item[\textit{ii)}] if $p=2$, then $G\cong D_8^{\circ k}$ or $G\cong Q_8\circ D_8^{\circ (k-1)}.$
\end{itemize}

If $G$ is an extraspecial $p$-group, then $o(x)\in\{ p, p^2\}$ for any non-trivial element $x$ of $G$. The results assembled in the following lemma can be found in \cite{13} (see Theorems 2.5 and 2.7).\\  

\textbf{Lemma 2.7.} \textit{Let $G$ be an extraspecial $p$-group. 
\begin{itemize}
\item[i)] If $G\cong M_{p^3}^{\circ k}$, then $n_p(G)=p^{2k}-1$ and $n_{p^2}(G)=(p-1)p^{2k}$;
\item[ii)] If $G\cong He_p^{\circ k}$, then $n_p(G)=p^{2k+1}-1$;
\item[iii)] If $G\cong D_8^{\circ k}$, then $n_2(G)=4^k+2^k-1$ and $n_4(G)=4^k-2^k;$
\item[iv)] If $G\cong Q_8\circ D_8^{\circ (k-1)}$, then $n_2(G)=4^k-2^k-1$ and $n_4(G)=4^k+2^k.$ 
\end{itemize}}

The last result of the paper shows the possible values of $|N_{n, L}|$ if $G$ is an extraspecial $p$-group. We also check whether $G\in\mathcal{D}$.\\

\textbf{Theorem 2.8.} \textit{Let $G$ be an extraspecial $p$-group, $n$ be a positive integer and $L\in\mathcal{P}^*(G)$.
\begin{itemize}
\item[a)] Assume that $p$ is odd.
\begin{itemize}
\item[i)] If $G\cong M_{p^3}^{\circ k}$, then $ |N_{n, L}|\in \{0, 1, 2^{p^{2k}}-1, (2^{p^{2k}}-1)^2, \ldots, (2^{p^{2k}}-1)^p,   2^{p^{2k+1}}-1\}$;
\item[ii)] If $G\cong He_p^{\circ k}$, then $|N_{n,L}|\in\{0, 1,  2^{p^{2k+1}}-1\}.$
\end{itemize}
Moreover, in both cases, $G\in\mathcal{D}$.
\item[b)] If $p=2$, then $|N_{n, L}|\in \{0, 1, 2^{4^k-2^k}-1, 2^{4^k+2^k}-1, (2^{4^k-2^k}-1)(2^{4^k+2^k}-1), 2^{2^{2k+1}}-1 \}$. Moreover, $G\not\in \mathcal{D}$.
\end{itemize}
}

\textbf{Proof.} \textit{a)} If $(n, exp(G))\in\{1, exp(G)\}$, then $|N_{n, L}|\in \{0, 1,  2^{p^{2k+1}}-1\}$ by Lemma 2.1 \textit{i)}, \textit{ii)}. Since $|\mathcal{P}^*(G)|^{|L|}=(2^{p^{2k+1}}-1)^{|L|}$, we also have $|N_{n,L}| \mid |\mathcal{P}^*(G)|^{|L|}$ if $|N_{n, L}|\neq 0$. This completely solves the case  $G\cong He_p^{\circ k}$ of exponent $p$.  

If $G\cong M_{p^3}^{\circ k}$, then $exp(G)=p^2$, so it remains to investigate the case $(n, exp(G))=p$. By reducing $n$ modulo $p^2$, we may assume that $n=mp$ with $m\in\{1, 2,\ldots, p-1\}$. Since $G'\subseteq Z(G)$, by Lemma 2.2 of \cite{3}, we know that $(xy)^i=x^iy^i[y,x]^{\frac{1}{2}i(i-1)}$ for any positive integer $i$ and any $x, y\in G$. By setting $i=n$ and taking into account that $p$ is odd and $G'\cong C_p$, we have that $$(xy)^n=x^ny^n, \forall \ x,y\in G.$$     
Hence, the function $f:G\longrightarrow G$ given by $f(x)=x^n, \ \forall \ x\in G$, is a group homomorphism. Assuming that $L\subseteq Im(f)$, by Lemma 2.1 \textit{iv)}, it suffices to evaluate $|Ker(f)|$. Clearly, $f$ is not trivial. On the other hand $G$ is extraspecial, so $x^p\in \Phi(G)$ for any $x\in G$. Then $\{ 1\}\lneq Im(f)\leq \Phi(G)\cong C_p$. It follows that $|Im(f)|=p$ and, consequently, $Ker(f)$ is a maximal subgroup of $G$. We conclude that if $(n, exp(G))=p$, then $|N_{n, L}|=(2^{p^{2k}}-1)^{|L|}$ with $|L|\in \{1, 2, \ldots, p\}.$ 
We also observe that $|N_{n, L}|$ divides $|\mathcal{P}^*(G)|^{|L|}$. This completely solves the case of extraspecial $p$-groups when $p$ is odd.

\textit{b)} If $(n, exp(G))\in\{1, exp(G)\}$, then $|N_{n, L}|\in \{0, 1,  2^{2^{2k+1}}-1\}$ by Lemma 2.1 \textit{i)}, \textit{ii)}. Since $exp(G)=4$, it remains to consider the case $(n, exp(G))=2$. We may assume that $n=2$. Let $l\in L$. If $o(l)=4$, then $|N_{n, l}|=0$. If $l=1$, according to Lemma 2.7 \textit{iii)}, \textit{iv)}, it is easy to deduce that
$$|N_{n,l}|=\begin{cases} 4^k+2^k & \text{, if } G\cong D_8^{\circ k}  \\  4^k-2^k & \text{, if } G\cong Q_8\circ D_8^{\circ (k-1)}\end{cases}.$$
Suppose that $o(l)=2$. Let $H=\langle x\rangle\cong C_4$ be a subgroup of $G$. Then, by Lemma 2.6 1., 2. of \cite{1}, it follows that $Z(G)\leq H$. Hence, if $\langle l\rangle=Z(G)$, then $l$ is the unique element of order 2 in any of the cyclic subgroups of order 4 of $G$. Consequently, if $o(x)=4$, then $x\in N_{n, l}$. According to Lemma 2.7 \textit{iii)}, \textit{iv)}, in this case, we have
$$|N_{n,l}|=\begin{cases} 4^k-2^k & \text{, if } G\cong D_8^{\circ k}  \\  4^k+2^k & \text{, if } G\cong Q_8\circ D_8^{\circ (k-1)}\end{cases}.$$
We also deduce that if $\langle l\rangle\neq Z(G)$, then $|N_{n, l}|=0$.
As a conclusion, if $(n, exp(G))=2$, then $|N_{n, L}|=0$ as long as $L$ contains at least one non-central element of $G$; otherwise, $L\in \mathcal{P}^*(Z(G))$ and $|N_{n, L}|\in \{ 2^{4^k-2^k}-1, 2^{4^k+2^k}-1, (2^{4^k-2^k}-1)(2^{4^k+2^k}-1)\}$. Thus, we  obtained all possible values of $|N_{n, L}|$. 

To show that $G\not\in \mathcal{D}$, consider $$L=\begin{cases} \{ 1\} & \text{, if } G\cong D_8^{\circ k}  \\  Z(G)\setminus \{1\} & \text{, if } G\cong Q_8\circ D_8^{\circ (k-1)}\end{cases}.$$
By our previous remarks, we have $|N_{n, L}|=2^{4^k+2^k}-1$ in both cases, while $|\mathcal{P}^*(G)|^{|L|}=2^{2^{2k+1}}-1$. Since $4^k+2^k \nmid 2^{2k+1}$, it follows that $|N_{n, L}|\nmid |\mathcal{P}^*(G)|^{|L|}$, so $G\not\in\mathcal{D}$.     
\hfill\rule{1,5mm}{1,5mm} 

\bigskip\noindent {\bf Acknowledgements.} The author is grateful to the reviewers for their remarks which improved the initial version of the paper.

\bigskip\noindent {\bf Declarations}

\bigskip\noindent {\bf Competing Interests.} The author has no relevant financial or non-financial interests to disclose. No data was collected in the course of this research.

\vspace*{3ex}
\small
\hfill
\begin{minipage}[t]{7cm}
Mihai-Silviu Lazorec \\
Faculty of  Mathematics \\
"Al.I. Cuza" University \\
Ia\c si, Romania \\
e-mail: {\tt silviu.lazorec@uaic.ro}
\end{minipage}

\begin{thebibliography}{100}
\bibitem{1} Bouc, S., Mazza, N., \textit{The Dade group of (almost) extraspecial p-groups}, J. Pure Appl. Algebra \textbf{192} (2004), 21-51.

\bibitem{2} Frobenius, G., \textit{\"Uber einen Fundamentalsatz der Gruppentheorie}, Berl. Sitz. (1903), 987-991.

\bibitem{3} Gorenstein, D., \textit{Finite groups}, 2nd edition, Chelsea Publishing Co., New York, 1980.

\bibitem{4} Hall, M., Jr., \textit{The theory of groups}, The Macmillan Company, New York, 1959. 

\bibitem{5} Hall, P., \textit{A Contribution to the Theory of Groups of Prime-Power Order}, Proc. London Math. Soc. (2) \textbf{36} (1934), 29-95.

\bibitem{6} Hall, P., \textit{On a Theorem of Frobenius}, Proc. London Math. Soc. (2) \textbf{40} (1935), no. 6, 468-501.

\bibitem{7} H\' ethelyi, L., L\' evai, L., \textit{On elements of order $p$ in powerful $p$-groups}, J. Algebra \textbf{270} (2003), no. 1, 1-6.

\bibitem{8} Huppert, B., \textit{Endliche Gruppen I}, Springer-Verlag, Berlin-New York, 1967.

\bibitem{9} Khurana, D., Khurana, A., \textit{A theorem of Frobenius and its applications}, Math. Mag. \textbf{78} (2005), no. 3, 220-225.

\bibitem{10} Klyachko, A. A., Mkrtchyan, A. A., \textit{How many tuples of group elements have a given property? With an appendix by Dmitrii V. Trushin}, Internat. J. Algebra Comput. \textbf{24} (2014), no. 4, 413-428.

\bibitem{11} Isaacs, I. M., Robinson, G. R., \textit{On a Theorem of Frobenius: Solutions of $x^n=1$ in Finite Groups}, Amer. Math. Monthly \textbf{99} (1992), no. 4, 352-354.

\bibitem{12} Iwasaki, S., \textit{A note on the nth roots ratio of a subgroup of a finite group}, J. Algebra \textbf{78} (1982), no. 2, 460-474.

\bibitem{13} Lazorec, M. S., \textit{Element orders in extraspecial groups}, Acta Math. Hung. \textbf{173} (2024), no. 2, 434-447.

\bibitem{14} Lazorec, M. S., T\u arn\u auceanu, M., \textit{Addendum to ``A generalization of a result on the sum of element orders of a finite group"}, Math. Slovaca \textbf{73} (2023), no. 1, 65-68.

\bibitem{15} Lubotzky, A., Mann, A., \textit{Powerful $p$-groups. I. Finite groups}, J. Algebra \textbf{105} (1987), no. 2, 506-515. 

\bibitem{16} Mazur, M., \textit{On powers in powerful $p$-groups}, J. Group Theory \textbf{10} (2007), no. 4, 431-433.

\bibitem{17} Sabatini, L., \textit{Products of subgroups, subnormality, and relative orders of elements}, Ars Math. Contemp. \textbf{24} (2024), no. 1, Paper No. 9, 9 pp.

\bibitem{18} Snapper, E., \textit{Normal subsets of finite groups}, Illinois J. Math. \textbf{13} (1969), 155-164.

\bibitem{19} T\u arn\u auceanu, M., \textit{A generalization of a result on the sum of element orders of a finite group}, Math. Slovaca \textbf{71} (2021), no. 3, 627-630.

\bibitem{20} The GAP Group, \textit{GAP -- Groups, Algorithms, and Programming}, Version 4.15.1; 2025. (\url{https://www.gap-system.org})
\end{thebibliography}
\end{document}